\documentclass[11pt]{amsart}
\usepackage[margin=1in]{geometry}
\usepackage{amsmath, amsfonts, amssymb, amsthm}
\usepackage{amsrefs}
\usepackage{mathrsfs}
\usepackage{enumitem} 
\usepackage{hyperref}
\usepackage{color}
\usepackage{dsfont}
\usepackage{csquotes}
\usepackage{epigraph}
\usepackage{comment}
\usepackage{mlmodern}
\usepackage[T1]{fontenc}

\usepackage{color}
\definecolor{darkblue}{rgb}{0.0,0.0,0.3}
\hypersetup{
    colorlinks=false,
    linkcolor=blue,
    urlcolor=darkblue,
    }

\newtheorem{theorem}{Theorem}[section]

\newtheorem{conjecture}[theorem]{Conjecture}

\newtheorem{lemma}[theorem]{Lemma}

\theoremstyle{definition}

\newtheorem{definition}[theorem]{Definition}

\theoremstyle{remark}
\newtheorem{remark}[theorem]{Remark}

\numberwithin{equation}{section}

\newcommand{\bR}{\mathbb{R}}
\newcommand{\bH}{\mathbb{H}}

\newcommand{\Garding}{G\r{a}rding}

\title[Hypersurfaces of constant curvature in hyperbolic space]{Curvature estimates for hypersurfaces of constant curvature in hyperbolic space II}
\author{Bin Wang}
\address[]{Department of Mathematics, The Chinese University of Hong Kong, Hong Kong, China.} 
\curraddr{Institute for Theoretical Sciences, Westlake University, Hangzhou, China.}
\email{wangbin@westlake.edu.cn}
\subjclass[2020]{Primary 53C21; Secondary 35J60, 53C40}
\keywords{Hessian equations, the asymptotic Plateau problem, hypersurfaces of constant curvature, fully nonlinear elliptic PDEs, a priori estimates.}

\begin{document}
\maketitle
\setcounter{tocdepth}{1} 
\begin{abstract}
In this note, we investigate the existence of smooth complete hypersurfaces in hyperbolic space with constant $(n-2)$-curvature and a prescribed asymptotic boundary at infinity. Previously, the existence was known only for a restricted range of curvature values, while in here, by deriving curvature estimates, we are able to deduce the existence for all possible curvature values.

\end{abstract}

\tableofcontents

\section{Introduction}
Fix $n \geq 3$. Let $\bH^{n+1}$ denote the hyperbolic space of dimension $n+1$ and let $\partial_{\infty} \bH^{n+1}$ denote its ideal boundary at infinity. Suppose $f: \bR^n \to \bR$ is a smooth symmetric function of $n$ variables satisfying some standard assumptions \cite[conditions (1.9)-(1.14) and condition (1.21)]{JEMS} on an open symmetric convex cone $K$ with vertex at the origin that contains the positive cone 
\[K^{+}:=\{\kappa \in \bR^n: \kappa_i>0 \quad \forall\ 1 \leq i \leq n\}.\]
Given a disjoint collection of closed embedded smooth $(n-1)$-dimensional submanifolds $\Gamma=\{\Gamma_1,\ldots,\Gamma_m\} \subseteq \partial_{\infty} \bH^{n+1}$ and a constant $0<\sigma<1$, we continue our previous investigations \cite{Bin-1,Bin-3} which follow the framework of Guan-Spruck \cite{JEMS} and seek a smooth complete hypersurface $\Sigma$ in $\bH^{n+1}$ satisfying 
\begin{equation}
\text{$f(\kappa[\Sigma])=\sigma$ and $\kappa[\Sigma] \in K$ on $\Sigma$}, \label{req1}
\end{equation}
with the asymptotic boundary
\begin{equation}
\partial \Sigma=\Gamma \label{req2}
\end{equation} at infinity; here $\kappa[\Sigma]=(\kappa_1,\ldots,\kappa_n)$ denotes the hyperbolic principal curvatures of $\Sigma$.

For the case when $K=K^{+}$, this problem has been extensively studied and eventually settled in a series of works \cite{JGA, JDG, SGAR} by Guan, Spruck, Szapiel and Xiao. Moreover, their results in \cite{JDG} are essentially optimal concerning that particular case. Here, we are concerned with the case when $K$ is a general cone for which the problem remains still not fully resolved. For this case, in a delicate work \cite{JEMS}, Guan-Spruck proved the existence of $\Sigma$ if $\sigma \in (\sigma_0,1)$, where $\sigma_0 \approx 0.37$. Later on, Xiao \cite{Xiao-B} improved the allowable range to $\sigma_0 \approx 0.14$ by refining their calculations. To this date, it is still not clear whether one can push $\sigma_0$ to be zero for a general curvature function $f$ in a general cone.

In this note, we choose the strategy that, instead of working for a full generality, we may first try some particularly interesting curvature functions $f$ in their natural cones $K$, with the hope that the attempt could be inspiring for future study. Specifically, we restrict ourselves to two classes of curvature functions $f=H_{k}^{1/k}$ or $f=(H_k/H_l)^{1/(k-l)}$, and our goal is to solve \eqref{req1}-\eqref{req2} for all $\sigma \in (0,1)$ in $K=K_k$, where
\[H_k(\kappa):=\binom{n}{k}^{-1}\sum_{1 \leq i_1<i_2<\cdots<i_k \leq n} \kappa_{i_1}\cdots \kappa_{i_k}\] is the normalized $k$-th elementary symmetric polynomial and 
\[K_k:=\{\kappa \in \bR^n: H_j(\kappa)>0 \quad \forall\ 1 \leq j \leq k\}\] is the $k$-th \Garding\ cone.

\begin{definition}
We say a smooth hypersurface $\Sigma$ in $\bH^{n+1}$ is strictly $k$-convex if $\kappa[\Sigma] \in K_k$ everywhere on the hypersurface.
\end{definition}

These two classes of curvature functions include a few notable examples and interpolate between them. For instance,
\begin{align*}
H_1(\kappa)&=\frac{1}{n}\sum_{i=1}^{n} \kappa_i \quad\quad\quad\quad\quad\quad\quad\phantom{i} \text{is the mean curvature},\\
H_2(\kappa)&=\frac{2}{n(n-1)}\sum_{i<j}\kappa_i\kappa_j \phantom{i}\quad\quad\quad \text{is the scalar curvature},\\
H_n(\kappa)&=\prod_{i=1}^{n}\kappa_i \quad\quad\quad\quad\quad\quad\quad\quad\phantom{i} \text{is the Gauss curvature}, \\
\left(\frac{H_n}{H_{n-1}}\right)(\kappa)&=n\left(\sum_{i=1}^{n} \frac{1}{\kappa_i}\right)^{-1} \quad\quad\quad\ \quad\phantom{i} \text{is the harmonic curvature};
\end{align*} other values of $k$ and $l$ also constitute remarkable geometric quantities; see \cite{Guan-Guan, Guan-Li-Li, Guan-Ma, Guan-Li, Sheng-Trudinger-Wang-CAG, Sheng-Trudinger-Wang-JDG}.

\begin{remark}
The quantity $H_k(\kappa)$ is called the $k$-th order mean curvature, or simply the $k$-curvature.
\end{remark}

The mean curvature case and the Gauss curvature case were resolved long ago in \cite{Anderson-1, Anderson-2, Lin, Tonegawa, Nelli-Spruck, AJM} and in \cite{Labourie, Rosenberg-Spruck}, respectively. In a recent work \cite{Bin-3}, we have settled the scalar curvature case and in an earlier investigation \cite{Bin-1}, we solved the problem for strictly $k$-convex hypersurfaces when $f=H_k/H_{k-1}$. In \cite{Lu}, Lu was able to obtain the existence of $\Sigma$ when $f=H_{n-1}^{1/(n-1)}$ and $K=K_{n-1}$. In this note, we claim that we can obtain the existence for $f=H_{n-2}^{1/(n-2)}$ in $K=K_{n-2}$.

\begin{theorem} \label{k=n-2}
Let $n \geq 5$.
Given a disjoint collection of closed embedded smooth $(n-1)$-dimensional submanifolds $\Gamma=\{\Gamma_1,\ldots,\Gamma_m\} \subseteq \partial_{\infty} \bH^{n+1}$ and a constant $0<\sigma<1$, if $\Gamma=\partial \Omega$ is mean-convex in the sense that its Euclidean mean curvature computed with respect to the inward normal vector is non-negative, then there exists a smooth complete strictly $(n-2)$-convex hypersurface $\Sigma$ in $\bH^{n+1}$ satisfying 
\[\text{$H_{n-2}^{1/(n-2)}(\kappa[\Sigma])=\sigma$ on $\Sigma$} \]
with the asymptotic boundary
\[\partial \Sigma=\Gamma \] at infinity.
\end{theorem}

To prove the theorem, we follow the framework of Guan-Spruck \cite{JEMS}, in which they reduced \eqref{req1}-\eqref{req2} to a Dirichlet boundary value problem for a fully nonlinear elliptic partial differential equation 
\begin{equation}
G(D^2u,Du,u)=\sigma, \quad u>0 \quad \text{in $\Omega \subseteq \bR^n$} \label{The PDE} 
\end{equation} 
with
\begin{equation}
u=0 \quad \text{on $\partial \Omega$}. \label{boundary condition}
\end{equation} As demonstrated in \cite{JEMS}, the PDE is degenerate when $u=0$ and so one must approximate the boundary condition by 
\begin{equation}
u=\varepsilon>0 \quad \text{on $\partial \Omega$}. \label{approximate boundary condition}
\end{equation}
The Dirichlet problem was solved for small $\varepsilon$ by proceeding with the standard continuity method for which one needs to establish a priori estimates for admissible solutions $u$ up to the second order. In fact, all estimates were perfectly established by Guan-Spruck \cite{JEMS} except for the global $C^2$ estimates, i.e.,
\begin{equation}
\max_{x \in \Omega} |D^2u| \leq C\left(1+\max_{\partial \Omega} |D^2u|\right), \label{global C2}
\end{equation} where $C>0$ must be independent of $\varepsilon$, which in their work was possible only for $\sigma \in (\sigma_0,1)$; this was the one and only place that they had to restrict the range of allowable $\sigma$'s. Thus, we shall focus on deriving \eqref{global C2}, or equivalently, on deriving a maximum principle for the hyperbolic principal curvatures $\kappa[\Sigma]$ for all $\sigma \in (0,1)$ and refer the reader to the original paper \cite{JEMS} for the complete theoretical framework.

In order to obtain \eqref{global C2} for achieving our Theorem \ref{k=n-2}, we will follow Lu's derivation \cite{Lu} in which he employed an inequality due to Ren-Wang \cite{Ren-Wang-1} that extracts enough positive terms from the concavity of the operator $f=H_{n-1}^{1/(n-1)}$ in $K=K_{n-1}$. These terms enabled Lu to deal with a crucial negative third order term when deriving \eqref{global C2}. The idea of exploring more concavity from the curvature function in concern may have come from Guan-Ren-Wang \cite{Guan-Ren-Wang}. Now, since Ren and Wang \cite{Ren-Wang-2} have just proved a similar concavity inequality for our operator $f=H_{n-2}^{1/(n-2)}$ in $K=K_{n-2}$, we will show that Lu's derivation can proceed in a similar fashion and the proof of Theorem \ref{k=n-2} follows.

\begin{remark}
In Theorem \ref{k=n-2}, the mean-convexity condition on $\Gamma$ was imposed to derive gradient estimates for admissible solutions; see \cite[Proposition 4.1]{JEMS}.
Also, we note that the dimensional constraint $n \geq 5$ is necessary for the concavity inequality to be effective, as demonstrated by the counterexample in \cite[p. 1297-1298]{Ren-Wang-1}. Hence, the theorem does not cover the case $n=4$ in which the curvature equation \eqref{req1} reduces to the constant scalar curvature equation in dimension four. The reader is referred to a recent work \cite{Bin-3} of ours for a treatment of the scalar curvature case in all dimensions.
\end{remark}

For the $k$-curvatures $f=H_{k}^{1/k}$ in $K=K_k$ with $3 \leq k \leq n-3$, the existence of solutions to \eqref{req1}-\eqref{req2} for a full range of $\sigma \in (0,1)$ remains open. In a recent preprint \cite{Hong-Zhang}, Hong and Zhang are able to derive \eqref{global C2} for constant $k$-curvature hypersurfaces that are semi-convex in the sense that $\kappa_i \geq -A$ everywhere for all $i$. Whether this additional assumption of semi-convexity can be removed is unclear. By defining $H_0:=1$, the curvature quotients $(H_k/H_l)^{1/(k-l)}$ can include the class of $k$-curvatures $f=H_{k}^{1/k}$ if we allow $l=0$. We state our ultimate goal as follows.

\begin{conjecture} \label{the conjecture}
Let $0 \leq l < k \leq n$.
Given a disjoint collection of closed embedded smooth $(n-1)$-dimensional submanifolds $\Gamma=\{\Gamma_1,\ldots,\Gamma_m\} \subseteq \partial_{\infty} \bH^{n+1}$ and a constant $0<\sigma<1$, if $\Gamma=\partial \Omega$ is mean-convex, then there exists a smooth complete strictly $k$-convex hypersurface $\Sigma$ in $\bH^{n+1}$ satisfying 
\[\text{$\left(\frac{H_k}{H_{l}}\right)^{\frac{1}{k-l}}(\kappa[\Sigma])=\sigma$ on $\Sigma$} \]
with the asymptotic boundary
\[\partial \Sigma=\Gamma \] at infinity.
\end{conjecture}

To conclude, we mention more related references for the interested reader. To solve the Dirichlet problem \eqref{The PDE}-\eqref{boundary condition}, there is a slightly different approximation scheme developed by Sui \cite{Sui-CVPDE}. Attempts that take this particular approach can be found in Sui-Sun \cite{Sui-CPAA, Sui-IMRN} and Jiao-Jiao \cite[Section 6]{Jiao-Jiao}. For the locally strictly convex case, i.e., when $K=K^{+}$, different approaches are also available. Using curvature flows, Xiao \cite{Xiao-A} provided a parabolic proof for the main results in \cite{JGA, SGAR}. In \cite{Hong-Li-Zhang}, Hong-Li-Zhang considered geodesic graphs whose boundaries rest upon equidistant hyperplanes instead of horospheres. A more geometric proof was recently provided by Smith \cite{Smith}. The proof applies a propagation of singularity result that is based on curvature estimates obtained by Sheng-Urbas-Wang \cite{Sheng-Urbas-Wang} and geometric properties of convex subsets of hyperbolic space.

The rest of this note is organized as follows. In section \ref{Preliminaries}, we list some basic formulas and fix notations. In section \ref{k=n-2 proof}, we derive the required curvature estimates and Theorem \ref{k=n-2} follows as a consequence of the theory in \cite{JEMS}. Theorem \ref{k=n-2} first appeared in an early preprint version of \cite{Bin-3}, it was then removed from \cite{Bin-3} and is entirely presented here instead.

\textbf{Errata to previous works.} In the first preprint version of this manuscript, we claimed that curvature estimates are possible for curvature quotients $(H_k/H_l)^{1/(k-l)}$ if $H_{k+1} \geq -A$. Now we wish to point out that there was a gap in the proof and hence it was removed from the current version. We would also like to point out that part (ii) and part (iii) should be removed from \cite[Corollary 1.2]{Bin-1}. The proof of curvature estimates in \cite[Section 4]{Bin-1} remains valid for corresponding cases. We sincerely apologize for these mistakes.

\section*{Acknowledgements}
The author would like to thank Professor Siyuan Lu for fruitful discussions on the subject and for suggesting the paper \cite{Ren-Wang-2} of Ren-Wang to us. This work was carried out while the author was still a Ph.D. student at The Chinese University of Hong Kong. The author is greatly indebted to Professor Man-Chun Lee and Professor Xiaolu Tan, for much support in various aspects academically.

\section{Preliminaries} \label{Preliminaries}
Let 
\[S_k:=\binom{n}{k} H_k\] denote the $k$-th elementary symmetric polynomial that is not normalized and we define
\[S_0=H_0:=1.\] Observe the property that
\[\frac{\partial S_k}{\partial \kappa_i}(\kappa)=S_{k-1}(\kappa_1,\ldots,\kappa_{i-1},0,\kappa_{i+1},\ldots,\kappa_n).\] We may denote the first and the second order derivatives by $S_{k-1}(\kappa|i)$ and $S_{k-2}(\kappa|pq)$, respectively. From the definition, the following can be readily deduced and verified.

\begin{lemma} \label{sigma-k properties 1}
For $\kappa = (\kappa_1,\ldots,\kappa_n) \in \bR^n$, we have
\begin{align*}
S_k(\kappa)&=\kappa_iS_{k-1}(\kappa|i)+S_k(\kappa|i),\\
\sum_{i=1}^{n} S_{k-1}(\kappa|i)&=(n-k+1)S_{k-1}(\kappa),\\
\sum_{i=1}^{n} \kappa_i S_{k-1}(\kappa|i)&=kS_k(\kappa).
\end{align*}
\end{lemma}

\begin{lemma} \label{sigma-k properties 2}
Let $1 \leq k \leq n$. The function $f(\kappa)=H_{k}^{1/k}(\kappa)$ satisfy the following in $K=K_k$: $f_i:=\partial f/\partial \kappa_i>0$ for all $1 \leq i \leq n$ and $f$ is concave. Moreover, $f>0$ in $K$ and $f=0$ on $\partial K$.
\end{lemma}

Following Guan-Spruck \cite{JEMS}, we use the half-space model for the hyperbolic space, 
\[\bH^{n+1}=\{(x,x_{n+1}) \in \bR^{n+1}: x_{n+1}>0\}\] that is equipped with the metric
\[ds^2=\frac{1}{x_{n+1}^2}\sum_{i=1}^{n} dx_{i}^2.\] Thus, the ideal boundary $\partial_{\infty} \bH^{n+1}$ is naturally identified with $\bR^n=\bR^n \times \{0\} \subseteq \bR^{n+1}$ and \eqref{req2} may be understood in the Euclidean sense. We will take $\Gamma=\partial \Omega$ for a smooth bounded domain $\Omega \subseteq \bR^n$.

Suppose $\Sigma=\{(x,u(x))\}$ can be represented as the graph of a function $u \in C^2(\Omega)$ over $\Omega$, which is oriented by the upward Euclidean unit normal vector field
\[\nu=\left(-\frac{Du}{w}, \frac{1}{w}\right), \quad \text{where $w=\sqrt{1+|Du|^2}$}.\] The component
\[\nu^{n+1}=\frac{1}{\sqrt{1+|Du|^2}}\] will prove to be useful in subsequent sections. The symbols $g$ and $\nabla$ will be used to denote the induced hyperbolic metric and the corresponding Levi-Civita connection on $\Sigma$. Since $\Sigma$ is also a submanifold of $\bR^{n+1}$, we will add a ``tilde'' over geometric quantities that are with respect to the Euclidean metric, $\tilde{g}$ and $\tilde{\nabla}$. In particular, some formulas \cite[Lemma 4.1]{Sui-CVPDE} that we will use is the following
\begin{gather}\label{geometric formula}
\begin{split}
\tilde{g}^{jk}&=\frac{\delta_{jk}}{u^2}, \quad h_{ij}=\frac{1}{u}\tilde{h}_{ij}+\frac{\nu^{n+1}}{u^2}\tilde{g}_{ij},\\
\nabla_i \nu^{n+1}&=-\tilde{h}_{ij}\tilde{g}^{jk}\nabla_k u, \quad \sum_{i=1}^{n} \frac{(\nabla_iu)^2}{u^2}=|\tilde{\nabla}u|^2=1-(\nu^{n+1})^2 \leq 1. 
\end{split}
\end{gather}

For an operator $F: \{\text{$n \times n$ symmetric matrices}\} \to \bR$ of the form
\[F(A)=f(\lambda(A)),\] where $\lambda(A)=(\lambda_1,\ldots,\lambda_n)$ are the eigenvalues of $A$ and $f: \bR^n \to \bR$ is a smooth symmetric function of $n$ variables. We denote
\[F^{ij}(A):=\frac{\partial F}{\partial a_{ij}}(A), \quad F^{ij,kl}(A):=\frac{\partial^2F}{\partial a_{ij} \partial a_{kl}}(A).\] In case the matrix $A=(a_{ij})$ is diagonal, we would have $F^{ij}=f_i\delta_{ij}$, where
\[f_{i}:=\frac{\partial f}{\partial \lambda_i}.\] Moreover, we have
\[\sum_{i,j} F^{ij}(A) a_{ij}=\sum_{i} f_{i}(\lambda(A))\lambda_i, \quad \sum_{i,j} F^{ij}(A)a_{ik}a_{jk}=\sum_{i} f_i(\lambda(A))\lambda_{i}^2.\] For the operator $F(A)=S_k(\lambda(A))$, we have
\[F^{ij}=S_{k-1}(\lambda|i)\delta_{ij}\] and 
\[F^{ij,pq}=\begin{cases}
S_{k-2}(\lambda|ip), & i=j, p=q, i \neq p \\
-S_{k-2}(\lambda|ip), & i=q, p=j, i \neq p \\
0, & \text{otherwise}
\end{cases}.\] It is common to abuse notations and use $S_k(A)$ to mean $S_{k}(\lambda(A))$, so we will also use $S_{k}^{ii}$ and $S_{k}^{pp,qq}$ to denote the first order and the second order derivatives, respectively.

\section{The curvature estimate} \label{k=n-2 proof}
In this section, we prove Theorem \ref{k=n-2}. As we have discussed in the Introduction, it is sufficient to derive, for admissible graphs, $\varepsilon$-independent curvature estimates that hold for all $\sigma \in (0,1)$; the rest will follow from the general framework of Guan-Spruck \cite{JEMS}. 

\begin{theorem}
Let $n \geq 5$ and $\sigma \in (0,1)$. Suppose $\Omega \subseteq \bR^n$ is a smooth bounded domain with a mean-convex boundary. For $f=H_{n-2}^{1/(n-2)}$ and $K=K_{n-2}$, let $u$ be a smooth solution of \eqref{The PDE} and \eqref{approximate boundary condition} whose graph $\Sigma$ is a smooth, strictly $(n-2)$-convex hypersurface in $\bH^{n+1}$. Then its largest principal curvature
\[\kappa_{\max}(X):=\max_{1 \leq i \leq n} \kappa_i(X)\] satisfies
\[\max_{\Sigma} \kappa_{\max} \leq C \left(1+ \max_{\partial \Sigma}\kappa_{\max}\right)\] for some universal constant $C>0$ that depends only on $n,\sigma$ and is independent of $\varepsilon$.
\end{theorem}

For the proof, let 
\[S_k:=\binom{n}{k}H_k,\] and we will work with the equation
\[S_{n-2}(\kappa)=\binom{n}{n-2}\sigma^{n-2}\] so that we can conveniently apply the following concavity inequality that is adapted from Ren-Wang's original version in \cite{Ren-Wang-2}.
\begin{lemma}
Let $F(A)=S_{n-2}(\kappa(A))$, where $A$ is an $n \times n$ symmetric matrix with eigenvalues $\kappa=(\kappa_1,\ldots,\kappa_n) \in K_{n-2}$ ordered as
\[\kappa_1 \geq \kappa_2 \geq \cdots \geq \kappa_n.\] Assume $n \geq 5$ and $\kappa_1$ is sufficiently large. If 
\[\kappa_1=\cdots=\kappa_m>\kappa_{m+1} \geq \cdots \geq \kappa_n \quad \text{for some $m \geq 1$},\]
then there exists some $C>0$ depending on known data, such that for an arbitrary vector $\xi=(\xi_1,\ldots,\xi_n) \in \bR^n$, we have that
\begin{gather}\label{Ren-Wang inequality}
\begin{split}
&-\sum_{p \neq q} \frac{F^{pp,qq}\xi_p \xi_q}{\kappa_1} - \frac{F^{11}\xi_{1}^2}{\kappa_{1}^2}+2\sum_{i > m}\frac{F^{ii}\xi_{i}^2}{\kappa_1(\kappa_1-\kappa_i)}\\ \geq &-\frac{C}{\kappa_1}\left(\sum_{i} F^{ii}\xi_{i}^2\right)^2-\frac{2}{\kappa_1}\sum_{1<i\leq m} F^{ii}\xi_{i}^2.
\end{split}
\end{gather}
When $m=1$, the index set $\{1<i\leq m\}$ is empty and the second sum on the right-hand side should be interpreted as zero.
\end{lemma}
\begin{remark}
By the vague phrase ``sufficiently large'', we mean that given any universal constant $C>0$ whose value depend only on the given data of the problem but not on $\varepsilon$, we are allowed to assume $\kappa_1 \geq C$.
\end{remark}
\begin{proof}[Proof of the lemma]
If $m=1$, then the proof is given in \cite[Lemma 2.11]{Bin-2}. We verify the inequality for $m>1$. Recall the original version from Ren-Wang \cite[Theorem 4]{Ren-Wang-2},
\[-\kappa_1\sum_{p,q} S_{n-2}^{pp,qq}\xi_p\xi_q - S_{n-2}^{11}\xi_{1}^2 + \sum_{i \neq 1} a_i \xi_{i}^2 \geq -C\kappa_1\left(\sum_i S_{n-2}^{ii}\xi_i\right)^2,\] where
\[a_i=S_{n-2}^{ii}+(\kappa_1+\kappa_i)S_{n-2}^{11,ii}.\] Then
\begin{align*}
&\ -\sum_{p \neq q} \frac{F^{pp,qq}\xi_p \xi_q}{\kappa_1} - \frac{F^{11}\xi_{1}^2}{\kappa_{1}^2}+2\sum_{i > m}\frac{F^{ii}\xi_{i}^2}{\kappa_1(\kappa_1-\kappa_i)}\\
\geq &\ -\frac{C}{\kappa_1}\left(\sum_i S_{n-2}^{ii}\xi_{i}\right)^2 - \frac{1}{\kappa_{1}^2}\sum_{i \neq 1} a_i\xi_{i}^2 + 2\sum_{i>m} \frac{F^{ii}\xi_{i}^2}{\kappa_1(\kappa_1-\kappa_i)}.
\end{align*}
For $1<i \leq m$, we have $\kappa_i=\kappa_1$. By the proof of Lemma 24 in \cite{Ren-Wang-2}, we have
\[a_i \leq 2\kappa_1 S_{n-2}^{ii} \quad \text{for $1 < i \leq m$}.\] While for $i>m$, we proceed as in \cite[Lemma 2.11]{Bin-2},
\begin{align*}
a_i&=S_{n-2}^{ii}+(\kappa_1+\kappa_i)S_{n-2}^{11,ii} \\
&=S_{n-2}^{ii}+\frac{\kappa_1+\kappa_i}{\kappa_1-\kappa_i} (S_{n-2}^{ii}-S_{n-2}^{11}) \\
&\leq \frac{2\kappa_1}{\kappa_1-\kappa_i}S_{n-2}^{ii} \quad \text{for $i>m$},
\end{align*} where we have noted that
\[\frac{\kappa_1+\kappa_i}{\kappa_1-\kappa_i}>0 \quad \text{for $n \geq 5$};\] see \eqref{negative kappa} below. Consequently, the lemma follows.
\end{proof}

\begin{proof}[Proof of the curvature estimate]
We consider the test function
\[Q= \frac{\kappa_{\max}}{(\nu^{n+1}-a)^{N}},\] where $N>0$ is some possibly large number to be determined later. Our arguments will mostly follow those of Lu in \cite{Lu}. Note, by the mean-convexity assumption on $\partial \Omega$ and the sharp gradient estimate in \cite[Proposition 4.1]{JEMS}, we have $\nu^{n+1} \geq \sigma$ everywhere on $\Sigma$.

Suppose $Q$ attains its maximum at some interior $X_0 \in \Sigma$. Let $\{\tau_1,\ldots,\tau_n\}$ be a local orthonormal frame around $X_0$ such that the second fundamental form $h_{ij}(X_0)=\kappa_i(X_0)\delta_{ij}$ is diagonalized and the principal curvatures are ordered as
\[\kappa_1 \geq \kappa_2 \geq \cdots \geq \kappa_n.\] If $Q$ is smooth at $X_0$, then we would have 
\begin{gather} \label{Q is smooth}
\begin{split}
\nabla_i\ \kappa_1&=\nabla_ih_{11},\\
\nabla_i \nabla_i\ \kappa_1&=\nabla_i\nabla_ih_{11}+2\sum_{p \neq 1} \frac{(\nabla_ih_{1p})^2}{\kappa_1-\kappa_p}.
\end{split}
\end{gather}
In case $\kappa_1$ has multiplicity $m>1$, i.e.,
\[\kappa_1=\kappa_2=\cdots=\kappa_m>\kappa_{m+1}\geq \cdots \geq \kappa_n,\] the function $Q$ would not be smooth. We can define a smooth function $\varphi(X)$ such that
\[Q(X_0)=\frac{\varphi(X)}{(\nu^{n+1}(X)-a)^N}.\] It follows that $\varphi \geq \kappa_1$ everywhere and $\varphi=\kappa_1$ at $X_0$. By a smooth approximation lemma \cite[Lemma 5]{BCD} due to Brendle-Choi-Daskalopoulos, we have
\begin{align}
\delta_{kl}\cdot \nabla_i\ \varphi &=\nabla_i h_{kl}, \quad 1 \leq k,l \leq m, \label{approximation}\\
\nabla_i\nabla_i\ \varphi&\geq \nabla_i\nabla_ih_{11}+2\sum_{p>m}\frac{(\nabla_ih_{1p})^2}{\kappa_1-\kappa_p}, \nonumber
\end{align} at the point $X_0$. Alternatively, one may apply a standard perturbation argument; see \cite{Chu,Lu}. From here onwards, all calculations will be done at the point $X_0$ without explicitly saying so, and we will use standard shorthand notations for the covariant derivatives, e.g., $h_{ijk}=\nabla_k h_{ij}$, $h_{ijkl}=\nabla_l\nabla_k h_{ij}$ and so on.
Now, if we were starting from the second case, then the critical equations will read as follows.
\begin{align}
0&=\frac{\varphi_i}{\varphi} - N\frac{\nabla_i\nu^{n+1}}{\nu^{n+1}} \nonumber \\
&= \frac{h_{11i}}{\kappa_1}-N\frac{\nabla_i \nu^{n+1}}{\nu^{n+1}}, \label{k=n-2 1st critical} \\
0&= \frac{\varphi_{ii}}{\varphi}-\frac{\varphi_{i}^2}{\varphi^2} - N\frac{\nabla_{ii}\nu^{n+1}}{\nu^{n+1}} + N \frac{(\nabla_{i}\nu^{n+1})^2}{(\nu^{n+1})^2} \nonumber\\
&\geq  \frac{h_{11ii}}{\kappa_1}+2\sum_{p>m}\frac{h_{1pi}^2}{\kappa_1(\kappa_1-\kappa_p)}-\frac{h_{11i}^2}{\kappa_{1}^2}- N\frac{\nabla_{ii}\nu^{n+1}}{\nu^{n+1}} + N \frac{(\nabla_{i}\nu^{n+1})^2}{(\nu^{n+1})^2}. \label{k=n-2 2nd critical 1}
\end{align} In either case, we would both get \eqref{k=n-2 1st critical} and \eqref{k=n-2 2nd critical 1}.

We may now assume $m>1$ and perform the derivation. The treatment for $m=1$ would be almost the same and simpler. Contracting \eqref{k=n-2 2nd critical 1} with $F^{ii}=S_{n-2}^{ii}$, we have
\begin{equation}
0 \geq  \sum_i \frac{F^{ii}h_{11ii}}{\kappa_1}+2\sum_i\sum_{p>m}\frac{F^{ii}h_{1pi}^2}{\kappa_1(\kappa_1-\kappa_p)}-\sum_i \frac{F^{ii}h_{11i}^2}{\kappa_{1}^2} - \frac{N}{\nu^{n+1}}\sum_i F^{ii}\nabla_{ii}\nu^{n+1}. \label{k=n-2 2nd critical 2}
\end{equation} Since $\bH^{n+1}$ has constant sectional curvature $-1$, by the Codazzi and Gauss equations, we have $h_{ijk}=h_{ikj}$ and 
\[h_{11ii}=h_{ii11}+\kappa_{1}^2\kappa_i-\kappa_1\kappa_{i}^2-\kappa_1+\kappa_i.\] Hence,
\[\sum_i F^{ii}h_{11ii}=\sum_i F^{ii}h_{ii11}-\kappa_1\left(\sum_i F^{ii}\kappa_{i}^2 + \sum_i F^{ii} \right) + (\kappa_{1}^2+1)(n-2)F,\] where we have also used Lemma \ref{sigma-k properties 1}. Furthermore,
by differentiating the curvature equation twice, we obtain that
\begin{equation}
\sum_i F^{ii}h_{ii1}=0 \label{k=n-2 differentiate once}
\end{equation} and
\begin{align*}
\sum_i F^{ii}h_{ii11}&=-\sum_{p,q,rs} F^{pq,rs}h_{pq1}h_{rs1} \\
&=-\sum_{p\neq q}F^{pp,qq}h_{pp1}h_{qq1}+\sum_{p \neq q}F^{pp,qq}h_{pq1}^2 \\
&\geq -\sum_{p\neq q}F^{pp,qq}h_{pp1}h_{qq1}+ 2 \sum_{i>m}\frac{F^{ii}-F^{11}}{\kappa_1-\kappa_i}h_{11i}^2.
\end{align*} Also, we have that
\begin{align*}
2\sum_i\sum_{p>m}\frac{F^{ii}h_{1pi}^2}{\kappa_1(\kappa_1-\kappa_p)} & \geq 2\sum_{p>m}\frac{F^{pp}h_{1pp}^2}{\kappa_1(\kappa_1-\kappa_p)}+2\sum_{p>m}\frac{F^{11}h_{1p1}^2}{\kappa_1(\kappa_1-\kappa_p)}\\
&=2\sum_{i>m}\frac{F^{ii}h_{ii1}^2}{\kappa_1(\kappa_1-\kappa_i)}+2\sum_{i>m}\frac{F^{11}h_{11i}^2}{\kappa_1(\kappa_1-\kappa_i)}.
\end{align*} Substituting all these back into \eqref{k=n-2 2nd critical 2}, we have
\begin{gather} \label{k=n-2 2nd critical 3}
\begin{split}
0&\geq -\sum_{p \neq q} \frac{F^{pp,qq}h_{pp1}h_{qq1}}{\kappa_1}+2\sum_{i>m}\frac{F^{ii}-F^{11}}{\kappa_1(\kappa_1-\kappa_i)}h_{11i}^2-\left(\sum_i F^{ii} + \sum_i F^{ii}\kappa_{i}^2\right)\\
&\quad +(n-2)F\kappa_1+2\sum_{i>m}\frac{F^{ii}h_{ii1}^2}{\kappa_1(\kappa_1-\kappa_i)}+2\sum_{i>m}\frac{F^{11}h_{11i}^2}{\kappa_1(\kappa_1-\kappa_i)}-\sum_i \frac{F^{ii}h_{11i}^2}{\kappa_{1}^2}\\
&\quad -\frac{N}{\nu^{n+1}}\sum_i F^{ii}\nabla_{ii} \nu^{n+1}.\\
\end{split}
\end{gather} 
For the third line in \eqref{k=n-2 2nd critical 3}, we proceed by a standard calculation \cite[Lemma 4.2]{Sui-CVPDE} and obtain that
\begin{align*}
&\ \sum_{i} F^{ii} \nabla_{ii}\nu^{n+1}\\
 =&\ 2\sum_i F^{ii}\frac{u_i}{u}\nabla_i\nu^{n+1}+(n-2)F[1+(\nu^{n+1})^2]\\
 & -\nu^{n+1}\left(\sum_i F^{ii}+\sum_i F^{ii}\kappa_{i}^2\right).
\end{align*} Therefore, from \eqref{k=n-2 2nd critical 3}, we get
\begin{gather} \label{k=n-2 2nd critical 4}
\begin{split}
0&\geq -\sum_{p \neq q} \frac{F^{pp,qq}h_{pp1}h_{qq1}}{\kappa_1} - \frac{F^{11}h_{111}^2}{\kappa_{1}^2}+2\sum_{i>m}\frac{F^{ii}h_{ii1}^2}{\kappa_1(\kappa_1-\kappa_i)} \\
&\quad +2\sum_{i>m}\frac{F^{ii}-F^{11}}{\kappa_1(\kappa_1-\kappa_i)}h_{11i}^2+2\sum_{i>m}\frac{F^{11}h_{11i}^2}{\kappa_1(\kappa_1-\kappa_i)}-\sum_{i\neq 1}\frac{F^{ii}h_{11i}^2}{\kappa_{1}^2}\\
&\quad +(N-1)\left(\sum_i F^{ii} + \sum_i F^{ii}\kappa_{i}^2\right)-2N\sum_i F^{ii}\frac{u_i}{u}\frac{\nabla_i\nu^{n+1}}{\nu^{n+1}}\\
&\quad +C(n,\sigma)\kappa_1-C(n)\cdot \sigma^{n-3} \cdot N,
\end{split}
\end{gather} where we have also made the following splitting,
\[\sum_{i} \frac{F^{ii}h_{11i}^2}{\kappa_{1}^2}=\frac{F^{11}h_{111}^2}{\kappa_{1}^2}+\sum_{i \neq 1} \frac{F^{ii}h_{11i}^2}{\kappa_{1}^2}.\]
According to Ren-Wang's concavity inequality \eqref{Ren-Wang inequality} and the one-time differentiation \eqref{k=n-2 differentiate once}, we have that
\begin{align*}
&-\sum_{p \neq q} \frac{F^{pp,qq}h_{pp1}h_{qq1}}{\kappa_1} - \frac{F^{11}h_{111}^2}{\kappa_{1}^2}+2\sum_{i>m}\frac{F^{ii}h_{ii1}^2}{\kappa_1(\kappa_1-\kappa_i)}\\ \geq &-\frac{C}{\kappa_1}\left(\sum_{i} F^{ii}h_{ii1}\right)^2-\frac{2}{\kappa_1}\sum_{1<i\leq m} F^{ii}h_{ii1}^2=0,
\end{align*} where we applied the concavity inequality to the vector $\xi_i=h_{ii1}$ and noted from \eqref{approximation} that
\[h_{ii1}=h_{i1i}=\delta_{i1} \cdot  (\varphi)_i = 0 \quad \text{for $1<i \leq m$}.\]
While for the second line in \eqref{k=n-2 2nd critical 4}, by again noting from \eqref{approximation}, we have
\begin{equation}
h_{11i}=h_{1i1}=\delta_{1i}\cdot (\varphi)_{1} = 0, \quad \text{for $1<i \leq m$}, \label{approximation 2}
\end{equation} and hence
\begin{align*}
&\ 2\sum_{i>m}\frac{F^{ii}-F^{11}}{\kappa_1(\kappa_1-\kappa_i)}h_{11i}^2+2\sum_{i>m}\frac{F^{11}h_{11i}^2}{\kappa_1(\kappa_1-\kappa_i)}-\sum_{i\neq 1}\frac{F^{ii}h_{11i}^2}{\kappa_{1}^2}\\
=&\ 2\sum_{i>m} \frac{F^{ii}h_{11i}^2}{\kappa_1(\kappa_1-\kappa_i)}-\sum_{i>m}\frac{F^{ii}h_{11i}^2}{\kappa_{1}^2}\\
=&\ \sum_{i>m}\left(\frac{2}{\kappa_1(\kappa_1-\kappa_i)}-\frac{1}{\kappa_{1}^2}\right) F^{ii}h_{11i}^2 \\
=&\ \sum_{i>m}\frac{\kappa_1+\kappa_i}{\kappa_1-\kappa_i} \frac{F^{ii}h_{11i}^2}{\kappa_{1}^2}\\
\geq &\ \sum_{i>m} \frac{n-4}{n}\frac{F^{ii}h_{11i}^2}{\kappa_{1}^2},
\end{align*} where we have used the following property \cite[Lemma 11]{Ren-Wang-2},
\begin{equation}
\kappa_i >-\frac{n-k}{k}\kappa_1 \quad \text{for $\kappa \in K_k$ with $\kappa_i \leq 0$}. \label{negative kappa}
\end{equation}
Note that $n-4 \geq 1$ by our assumption $n \geq 5$. Thus, \eqref{k=n-2 2nd critical 4} becomes
\begin{gather}\label{k=n-2 2nd critical 5}
\begin{split}
0 & \geq (N-1)\left(\sum_i F^{ii} + \sum_i F^{ii}\kappa_{i}^2\right)+C(n,\sigma)\kappa_1-C(n)\sigma^{n-3}N\\
&\quad + \frac{1}{n}\sum_{i>m} F^{ii}\frac{h_{11i}^2}{\kappa_{1}^2} - 2 N\sum_i F^{ii}\frac{u_i}{u}\frac{\nabla_i \nu^{n+1}}{\nu^{n+1}}.
\end{split}
\end{gather} We continue to deal with the second line in \eqref{k=n-2 2nd critical 5}. By the first order critical equation \eqref{k=n-2 1st critical}, when $\kappa_1(X_0)>1 \geq \nu^{n+1}(X_0)$, we have that
\begin{align*}
&\ \frac{1}{n}\sum_{i>m} F^{ii}\frac{h_{11i}^2}{\kappa_{1}^2} - 2 N\sum_i F^{ii}\frac{u_i}{u}\frac{\nabla_i \nu^{n+1}}{\nu^{n+1}} \\
=&\ \frac{N^2}{n}\sum_{i>m} F^{ii}\left(\frac{\nabla_i\nu^{n+1}}{\nu^{n+1}}\right)^2-2 N\sum_i F^{ii}\frac{u_i}{u}\frac{\nabla_i \nu^{n+1}}{\nu^{n+1}} \\
=&\ \frac{N^2}{n}\sum_{i>m} F^{ii}\frac{u_{i}^2}{u^2}\left(\frac{\kappa_i-\nu^{n+1}}{\nu^{n+1}}\right)^2+2 N\sum_{i} F^{ii}\frac{u_{i}^2}{u^2}\frac{\kappa_i-\nu^{n+1}}{\nu^{n+1}} \\
\geq &\ \frac{N^2}{n}\sum_{i>m} F^{ii}\frac{u_{i}^2}{u^2}\left(\frac{\kappa_i-\nu^{n+1}}{\nu^{n+1}}\right)^2+2 N\sum_{i>m} F^{ii}\frac{u_{i}^2}{u^2}\frac{\kappa_i-\nu^{n+1}}{\nu^{n+1}} \\
> &\ -2n\sum_{i>m} F^{ii},
\end{align*} where we have used \eqref{geometric formula} to get the following
\[\nabla_i\nu^{n+1}=\frac{u_{i}}{u}(\nu^{n+1}-\kappa_i), \quad \frac{u_{i}^2}{u^2} \leq \sum_i \frac{u_{i}^2}{u^2} \leq 1.\]

Therefore, by choosing $N$ large enough, e.g., $N=2n+1$, we have that
\[(N-1)\sum_i F^{ii}+\frac{1}{n}\sum_{i>m} F^{ii}\frac{h_{11i}^2}{\kappa_{1}^2} - 2 N\sum_i F^{ii}\frac{u_i}{u}\frac{\nabla_i \nu^{n+1}}{\nu^{n+1}} \geq 0.\] Hence, from \eqref{k=n-2 2nd critical 5}, we are left with
\[0 \geq C(n,\sigma) \kappa_1-C(n)\sigma^{n-3}N\] and the desired estimate follows.
\end{proof}

\bibliography{refs}
\end{document}